\documentclass{article}
\usepackage[utf8]{inputenc}

\usepackage{latexsym} 
\usepackage{amssymb} 
\usepackage{amsmath}
\usepackage{amsthm}
\usepackage{enumerate} 
\usepackage{graphics}
\usepackage{graphicx}

\newtheorem{thm}{Theorem}

\newtheorem{example}[thm]{Example}

\newcommand{\pf}{\noindent{\bf Proof}\enskip}

\newcommand{\F}{{\mathbb F}}

\newcommand{\Q}{{\cal Q}}

\title{\vskip-3cm Covering all but the low weight vertices of the unit cube}
\author{P. Sziklai and Zs. Weiner}
\date{}

\begin{document}

\maketitle

\begin{abstract}
In this paper we discuss a result similar to the polynomial version of the Alon-Füredi theorem \cite{AF}. We prove that if you want to cover the vertices of the $n$-dimensional unit cube, except those of weight at most $r$ then you need an algebraic surface of degree at least $n-r$.    
\end{abstract}


\section{Introduction}

Let $\Q$ be the unit cube $\{0,1\}^n$ of the vector space $\F^n$, where $\F$ is a field. There is a quadratic surface covering all the vertices of $\Q$. But if we forbid to cover some of the vertices it becomes a much more difficult question how (i.e. by how small degree polynomial) can we achieve it. A typical result of this flavour states that if we forbid one vertex (e.g. the origin) then we need a polynomial of degree at least $n$; or more generally, formulated the other way around in \cite{AF}, if a polynomial of degree $d$ does not vanish completely on the grid $S_1\times...\times S_n$, where $0<|S_i|, S_i\subset\F\ \forall i$, then it is nonzero on at least $\min\prod y_i$ points of the grid, where the minimum is taken over all sets of integers $0<y_i\leq |S_i|\ \forall i$, the sum of which is at least $\sum |S_i| -d$.

There is an abundance of results related to the Alon-Füredi paper, we do not survey them here.

\section{The main result}

The {\it weight} of a vector is just the number of nonzero coordinates of it. The next theorem extends the result of Alon-Füredi \cite{AF}.

\begin{thm}
In ${\F}^n$, if for a polynomial $f\in\F[x_1,x_2,...,x_n]$ of degree $d$, we have $f(x)=0$ for each vertex $x$ of the unit cube except the vertices of weight $\le r$, where $f(x)\ne 0$, then $d\geq n-r$.
\end{thm}
Note that the theorem is sharp, an obvious example is the following polynomial  (and there are many others).

\begin{example}
If $\mathrm{char}\ \F=0$ or $n<\mathrm{char}\ \F$ then\\
$\displaystyle{f(x_1,x_2,...,x_n)=\prod_{s=r+1}^n (x_1+x_2+...+x_n-s)}$ is a polynomial vanishing on the vertices of the unit cube of weight at least $r+1$ and nonzero on the rest.
\end{example}

There are many versions and proofs of similar results, see \cite{AF}. Here we show one,  which depends on careful examination of the coefficients of the polynomial.

\medskip
\pf of the theorem.
Suppose that, on the contrary, $d<n-r$.
Write $$f(x_1,x_2,...,x_n)=\sum_{0\leq i_1+i_2+...+i_n\leq d} a_{i_1,i_2,...,i_n} x_1^{i_1}x_2^{i_2}...x_n^{i_n}\ \ .$$
We say that a term {\it contains} the variable $x_k$ if the exponent of $x_k$ in the term is nonzero. 
Define $\alpha_{\{j_1,j_2,...,j_s\}}$ or $\alpha_{j_1,j_2,...,j_s}$ as the sum of the coefficients of the terms of $f$, {\it containing} precisely the variables $x_{j_1}, x_{j_2}, ..., x_{j_s}$ (i.e. with exponent at least 1) but no other variables. Note that our assumption $d<n-r$ implies that
$$\alpha_J=0 \textrm{ for all } J\subset\{1,...,n\},\ |J|\geq n-r. \eqno{(*)}$$

Substituting vertices of $\Q$ with weight $\leq r$ (i.e. vectors with at most $r$ coordinates being 1 and all the others zero), we get that 
$$\alpha_J\ne -\sum_{A\subsetneq J} \alpha_A \ \ \ \ \ \ \textrm{ for } 1\leq s\leq r,\ J\subseteq \{1,...,n\},\ |J|=s\ \ .$$

Now substituting vertices of $\Q$ with weight $s$, where $0\leq s \leq n$, and denoting $r^*=\min(s,r)$, by Möbius-inversion we get that for $J\subseteq \{1,...,n\},\ |J|=s$ 
$$\alpha_J = \sum_{A\subseteq J} (-1)^{|J\setminus A|} f(A) = \sum_{u=0}^{r^*} (-1)^{s-u} \sum_{A\subseteq J\atop |A|=u} \sum_{B\subseteq A}\alpha_B =$$ 
$$\sum_{u=0}^{r^*} (-1)^{s-u} \sum_{B\subseteq J\atop (|B|\leq u)} {s-|B|\choose u-|B|}\alpha_B=
\sum_{t=0}^{r^*} \Bigg(\sum_{u=t}^{r^*} (-1)^{s-u} {s-t\choose u-t}\Bigg) \sum_{B\subseteq J\atop |B|=t} \alpha_B\ . \eqno{(**)}$$
As\\  
$$\sum_{u=t}^{r^*} (-1)^{s-u} {s-t\choose u-t} =\left\{ \begin{array}{lr}
1 & \text{ if } t=s=r^*;\\
0 & \text{ if } 0\leq t<s=r^*; \text{ and}\\
(-1)^{s-r^*} {s-1-t\choose r^*-t} & \text{ otherwise;}
\end{array} \right.
$$ 
from $(**)$ we have in the case $s\leq r$ (the obvious)\\
$$\alpha_J\ =\ \alpha_J\ \ ;$$
while in the case $r<s\leq n$ we get
$$\alpha_J\ =\ \sum_{t=0}^r (-1)^{s-r} {s-1-t\choose r-t} \sum_{B\subseteq J\atop |B|=t} \alpha_B\ \ .$$
This is a set of linear equations, and its equations can be indexed by the complement sets $\bar{J}=\{1,...,n\}\setminus J$ and the "variables" are the coefficient sums $\alpha_B$ for the subsets $B\subseteq\{1,...,n\},\ \ |B|\leq r$. If we consider the equations $|\bar{J}|\leq r$ then we get a system of {\it homogeneous} linear equations of size $\sum_{i=0}^r{n\choose i}\  \times\ \sum_{i=0}^r{n\choose i}$, as the corresponding $\alpha_J$ values on the {\it left hand sides} are all zero by ($*$). \\

{\bf Firstly}, suppose that $r<n/2$.

The rows and the columns of the matrix $M$ of this system of equations are indexed by the subsets of size at most $r$ of $\{1,...,n\}$, and an entry $m_{A,B}$ is equal to $(-1)^{n-r-|A|}{n-1-|A|-|B| \choose r-|B|}$ whenever $A$ and $B$ are disjoint subsets, and zero otherwise. \\

\noindent{\bf Claim:} $M=M^{-1}$.

Proof: in $MM$, the entry indexed by the subsets $A$ and $B$ is the following:

if $A=B$ then 
$$\sum_U m_{A,U}m_{U,A}=
\sum_{U\subseteq\bar{A}}(-1)^{|A|+|U|}{n-1-|A|-|U| \choose r-|U|}{n-1-|A|-|U| \choose r-|A|}=$$
\begin{equation}\label{zeilberger1}
(-1)^{|A|}\sum_{u=0}^{\min(n-|A|,r)} (-1)^u {n-|A|\choose u}{n-1-|A|-u \choose r-u}{n-1-|A|-u \choose r-|A|}= 1\ .\end{equation}

If $A\ne B$ then 
$$\sum_U m_{A,U}m_{U,B}=
\sum_{U\subseteq\overline{A\cup B}}(-1)^{|A|+|U|}{n-1-|A|-|U| \choose r-|U|}{n-1-|B|-|U| \choose r-|B|}=$$
 \begin{equation}\label{zeilberger2}(-1)^{|A|}\sum_{u=0}^{\min(n-|A\cup B|,r)} (-1)^u {n-|A\cup B|\choose u}{n-1-|A|-u \choose r-u}{n-1-|B|-u \choose r-|B|}= 0\ .\end{equation}

These equalities can be proved by Zeilberger's method (see the Appendix), we used the fastZeil Mathematica package developed by Paule, Schorn and Riese \cite{Paule}. We are grateful for them to share the package with us and for their helpful advice.

Hence $M$ is invertible indeed and the unique solution is $\alpha_J=0$ for all $|J|\leq r$. But this is a contradiction.\\

{\bf Secondly}, suppose that $r\geq n/2$.

Now the matrix $M$ is similar, but (as we have now equations for $n-r\leq s\leq r$),
it contains rows belonging to equations $\alpha_J=\alpha_J$, i.e. 
in the row indexed by $A=\bar{J}$, $|J|=s,\ n-r\leq |A|\leq r$, the element
$m_{A,B}=1$ for $B=\bar{A}$ and zero otherwise.

The rows and the columns of the matrix $M$ of this system of equations are still indexed by the subsets of size at most $r$ of $\{1,...,n\}$, and the rows indexed by sets of size less than $n-r$ remained the same, i.e. the entry $m_{A,B}$ is equal to $(-1)^{n-r-|A|}{n-1-|A|-|B| \choose r-|B|}$ whenever $A$ and $B$ are disjoint subsets, and zero otherwise. \\
Note that if we order the index sets increasingly w.r.t. their size, and in the same way for rows and columns, then in $M$ we can see an $\sum_{i=n-r}^r {n\choose i}\times \sum_{i=n-r}^r {n\choose i}$ identity matrix in the bottom-right corner, only zeroes on its left, and in the upper-left corner we find $M_0$ of size $\sum_{i=0}^{n-r-1} {n\choose i}\times \sum_{i=0}^{n-r-1} {n\choose i}$ which is similar to the 'old' version of $M$ above and we can prove $M_0=M_0^{-1}$.

It follows that $M$ is invertible indeed and the unique solution is $\alpha_J=0$ for all $|J|\leq r$. But this is a contradiction again.
\qed
\\

We note that in the extremal case $d=n-r$ the same equalities can be used to describe the $\alpha_J$-s; there remains a lot of freedom to choose the coefficients of $f$.

\section{Appendix}
Here we sketch the proof of the two equalities (\ref{zeilberger1}) and (\ref{zeilberger2}) which serve the proof of $M = M^{-1}$.
Note that for $r=0$, the matrix $M$ is 1-by-1 with its only entry being $(-1)^n$; while for $r=1$ we have an $(n+1)\times(n+1)$ matrix for which, again, it is easy to check (\ref{zeilberger1}) and (\ref{zeilberger2}).

Now to prove (\ref{zeilberger1}) let
$$S_1(r) = \sum_{u=0}^r(-1)^{u+|A|} {n-|A|\choose u}{n-1-|A|-u \choose r-u}{n-1-|A|-u \choose r-|A|}.$$ Note that in (\ref{zeilberger1}) the sum runs until $\min(n-|A|,r)$ which is $r$
as $r < n/2$.
We want to show that $S_1(r)=1$, for $r<n/2$. Let $n-|A|=m$ and $|A| = a$.
Zeilberger's method provides the recursion:
\begin{multline*} 
 -(a - r - 1) (m - r - 1)(a + m - 2 r - 4) (a + m - r - 1)\ S_1(r)+ \\ (a + m - 2 r - 3) (a^2 m - a^2 r - a^2 + a m^2 - 2 a m r -2 a m + a r^2 + a r - a - m^2 r - m^2 + \\ m r^2 + m r - m + 2 r^2 + 6 r + 4)\ S_1(r+1)- \\- (r + 2) (a - r - 2) (m - r - 2) (a + m - 2 r - 2)S_1(r+2) =0
\end{multline*}
For $r,a < n/2$, the coefficient of $S_1(r+2)$ is nonzero. From the first paragraph of this section, $S_1(r) = 1$ for $r=0, 1$ and so, comparing the coefficients of $S_1(r), S_1(r+1)$ and $S_1(r+2)$ we get, by induction, that $S_1(r)=1$ for all $r$.

In order to prove (\ref{zeilberger2}) let
$$S_2(r) = \sum_{u=0}^r(-1)^{u+|A|} {n-|A\cup B|\choose u}{n-1-|A|-u \choose r-u}{n-1-|B|-u \choose r-|B|}.$$ In (\ref{zeilberger2}), the sum runs until $\min(n-|A\cup B|,r)$, but when $u > n-|A\cup B|$ then ${n-|A\cup B|\choose u}=0$, so the result does not change if we sum up to $r$.
We want to show that $S_2(r)=0$, for $r<n/2$. Let $n-|A\cup B| = m$, $|A\cap B| = w$, $|A| = a$ and  $|B| = b$.
Zeilberger's method provides the recursion: 
\begin{multline*} 
 -(a - r - 1)(b + m - r - w - 1) (a + b + m - 2 r - w - 4) (a + b + m - r - w - 1) S_2(r)\\  -(a + b + m - 2 r - w - 3)
(a^2 b - a^2 r - a^2 w - 2 a^2 + a b^2 + a b m - 2 a b r - 3 a b w - 4 a b - 2 a m w \\- a m + a r^2 +
4 a r w + 5 a r + 2 a w^2 + 7 a w + 5 a - b^2 r - b^2 w - 2 b^2 - 2 b m w - b m + b r^2 +
4 b r w \\+ 5 b r + 2 b w^2 + 7 b w + 5 b + m^2 r - m^2 w + m^2 - m r^2 + 2 m r w - m r +
2 m w^2 + 4 m w + m \\- 3 r^2 w - 2 r^2 - 3 r w^2 - 11 r w - 6 r - w^3 - 5 w^2 - 9 w - 4)S_2(r+1) \\+
(r + 2) (b - r - 2) (-a - m + r + w + 2) (a + b + m - 2 r - w - 2)S_2(r + 2) = 0.
\end{multline*}

Again we see that $S_2(r)=0$ for $r=0, 1$ and the coefficient of $S_2(r+2)$ is nonzero when $r<n/2$ and so $S_2(r)$ is always $0$.

\section{Addendum}

After publication of this paper, the authors learned that a more general version of their result had been proved independently, slightly earlier, by Venkitesh \cite{Venkitesh}, Corollary 33. In \cite{Venkitesh}, this is a corollary of a nice, rather complex series of results, so our 2 or 3 pages long proof remains still interesting; and we believe that this application of Zeilberger's method is still worth publishing.

\section{Acknowledgements}

The second author acknowledges the partial support of the National Research, Development and Innovation Office – NKFIH, grant no. K 124950. The first author is grateful for the partial support of project K 120154 of the National Research, Development and Innovation Fund of Hungary; and for the support of the National Research, Development and Innovation Office within the framework of the Thematic Excellence Program 2021 - National Research Subprogramme: “Artificial intelligence, large networks, data security: mathematical foundation and applications”.

\bigskip

Peter Sziklai 

ELTE  Eötvös Loránd University, Budapest, Hungary

\texttt{peter.sziklai@ttk.elte.hu}

\medskip

Zsuzsa Weiner

ELKH-ELTE GAC Research Group, Budapest, Hungary

\texttt{zsuzsa.weiner@gmail.com}


\begin{thebibliography}{}
\bibitem{AF}
N. Alon and Z. Füredi, Covering the cube by affine hyperplanes, {\it European J. Combinatorics}, 14:79–83, 1993.

\bibitem{Paule}
P. Paule and M. Schorn, A Mathematica Version of Zeilberger’s
Algorithm for Proving Binomial Coefficient Identities,
{\it J. Symbolic Comput.}, 20:673-698, 1995.

\bibitem{Venkitesh}
S. Venkitesh, Covering Symmetric Sets of the Boolean Cube by Affine Hyperplanes,
{\it Electronic J. Combin.} Vol. 29, Issue 2 (2022), Paper 2.22 (31 pages)

\end{thebibliography}
\end{document}